\documentclass{article}
\usepackage{amsfonts}
\usepackage{amssymb}
\usepackage{amsmath}
\usepackage{amsthm}
\usepackage{tikz-cd}
\usepackage{float}
\usepackage{tikz}
\usetikzlibrary{trees}
\usepackage{pifont}
\usepackage{enumitem}

\theoremstyle{plain}
\newtheorem{thm}{Theorem}

\newtheorem{defn}[thm]{Definition}
\theoremstyle{definition}
\newtheorem{ej}{Example}
\newtheorem*{af}{Claim}

\newtheorem{coro}[thm]{Corollary}
\theoremstyle{remark}
\newtheorem*{nota}{Remark}

\newcommand{\ga}{\alpha}

\newcommand{\gd}{\delta}
\newcommand{\gl}{\lambda}
\newcommand{\gk}{\kappa}
\newcommand{\gf}{\varphi}
\newcommand{\gS}{\Sigma}

\renewcommand{\gg}{\gamma}

\newcommand{\go}{\omega}

\newcommand{\e}{\exists\,}

\newcommand{\po}{{\mathcal L_{\omega\omega}}}
\usepackage{xspace}

\newcommand{\suc}[1]{\left\langle#1\right\rangle}

\newcommand{\vy}{\vec{y}}

\newcommand{\vx}{\vec{x}}

\newcommand{\uhr}{\upharpoonright}
\newcommand{\vfi}{\varphi}
\newcommand{\lra}{\leftrightarrow}

\newcommand{\LA}{\mathcal L}

\newcommand{\A}{\mathfrak A}
\newcommand{\B}{\mathfrak B}

\newcommand{\ssi}{\Longleftrightarrow}

\newcommand{\mf}{\mathfrak }
\newcommand{\mrm}{\mathrm }

\newcommand{\mtx}[1]{\begin{matrix}#1 \end{matrix}}

\newcommand{\qr}{\mrm{qr}}

\newenvironment{prueba}{\paragraph{\emph{Claim's Proof.} 
\normalfont}}{\hfill$\blacksquare$}

\begin{document}

\title{Cardinal Transfer Events for Infinitary Logic with Infinitary Quantification}
\author{Edgar A. Valenzuela Nuncio}


\maketitle

\begin{abstract}
Assuming the existence of a certain theory $\Phi$, we solve the Gap-$n$ problem for more than one predicate. This will be done for infinitary logic with infinitary quantification. 
\end{abstract}

\section{Introduction}

The Gap-$1$ problem states that, given a vocabulary $\tau$ which includes a predicate $U$, if we can find a model $\A$ with power $\gk^+$ and such that $|U^\A|=\gk$ then we can find an elementary equivalent $\B$ model of power $\nu^+$ such that $|U^\B|=\nu$. We say that the type of $\A$ is $(\gk^+,\gk)$ and $(\nu^+,\nu)$ is the type of $\B$. With
$$
(\gl,\gk) \to (\nu,\mu),
$$
we signifiy that the existence of a model $\A$ of type $(\gl,\gk)$ implies the existence of a model $\B\equiv \A$ of type $(\nu,\mu)$.
In general, the Gap-$n$ problem states that
$$
(\gk^{(n)},\mu) \to (\nu^{(n)},\nu).
$$
 In the case $n=1$ there are positive answers using $L$ or $L$-like hypotheses, such as the GCH \cite{Chang2} or \emph{Morasses} \cite{Jen}, \cite{Vi06}, \cite{Vi10}, \cite{Vi}. Also, in the case $n=2$ the solutions make use of $L$-like hypotheses, \cite{Dev}, \cite{Vill}. The route we follow goes in the opposite direction. We will assume the existence of a theory $\Phi$ that already implies $V\neq L$.

 Another way to generalize this notion is to include several predicates $U_0,...,U_k$. In this case, one can study cardinal transfer problems such as
 $$
(\gk^{(n_k)},...,\gk^{(n_0)})\to (\gl^{(n_k)},...,\gl^{(n_0)}),
 $$
 
  for $n_0<...<n_k$. This generalization will allow us to study the behavior of the possible sizes when using several predicates. 
  
  Also, we could consider this transfer theorems for extensions of first-order logic. And we will do so in this work. An infinitary version of the Gap transfer problem will be proved here making use of the existence of a certain theory $\Phi$.
  
   The infinitary logic, $\LA'_{\gk\gl}$, we will be dealing with depends on the rank of quantification for formulas. Common uses of infinitary logic have the restriction $\gl\leq\gk$. The logic used here does not have this restriction.


In the first section of this paper we will present the logic $\LA'_{\gk\gl}$. Then there will be shown some relevant results around the Gap-$n$ problem. The results obtained will be for the levels $L_\gl$ of the constructible hierarchy. We will see that, assuming the existence of $0^\sharp$, for a vocabulary $\tau$ which includes the predicates $U_0,....,U_{n-1}$,  $\mu_0\leq ...\leq \mu_{n-1}$ and $\nu_0\leq...\leq \nu_{n-1}$ we obtain
$$
(\mu_{n-1},...,\mu_0)\to(\nu_{n-1},...,\nu_0).
$$

Particularly, we will obtain
$$
(\mu_{k-1}^{(n_k)},...,\mu_0^{(n_0)})\to(\nu_{k-1}^{(n_k)},...,\nu_0^{(n_0)}),
$$
for many instances. 

The use of the models $L_\gl$ in the first section is by no means trivial. We will use them because the results of the previous section will allow us to generalize the Gap-$n$ problem in the second section of this work. We will define a generalization of the Gap-$n$ problem in order to consider several predicates which we name \emph{Gap-$(n_0,...,n_k)$ problem}. Assuming the existence of an adequate theory $\Phi$, we will be able to obtain the generalization of the previous section but now for a larger class of vocabularies. This is the content of Theorem \ref{main}. The importance of the theory $\Phi$ is that it will generate a set of indiscernibles for the constructible relativized universe $L[a]$, for some set $a\subset L$.

 This will yield the upwards version of the Gap-$(n_0,...,n_k)$ problem for infinitary logic. Particularly, for first-order logic, we will obtain the full Gap-$(n_0,...,n_k)$ outside of the constructible universe $L$.

\section{The Logic $\LA'_{\gk\gl}$}

We begin by defining the  relevant concepts starting with the logic $\LA_{\gk\gl}'$. To define this logic first we need to clarify what will be the difference between $\LA_{\gk\gl}'$ and $\LA_{\gk\gl}$. In order to do so, we need to work with the \emph{quantifier rank} of a formula. 

\begin{defn}
Fix any $\vfi\in \LA_{\infty\infty}$, the \emph{quantifier rank of} $\vfi$, written as  $\qr(\gf)$, is defined recursively as shown next.
$$
\mtx{
\qr(\gf) = 0 & \gf \text{ atomic,}\\
\qr(\neg\gf) = \qr(\gf), &  \\
\qr(\bigwedge \Phi) = \sup_{\gf\in \Phi}\qr(\gf), & \\
\qr(\forall \vx \gf(\vx)) = \qr(\e \vx\gf(\vx)) &  = \qr(\gf)+1.
}
$$
\end{defn}

With this at hand, we are able to define our logic.

\begin{defn}
Let $\LA'_{\lambda\gk}\leq \LA_{\infty\infty}$ be the set of formulas with $<\kappa$ free variables and with quantifier rank $\leq\lambda$.
\end{defn}

The next examples show the difference between $\LA_{\gk\gl}$ and $\LA'_{\gk\gl}$.

\begin{ej}\ 
\begin{enumerate}

\item Notice that the index $\gk$ of $\LA'_{\gl\gk}$ deals with the amount of conjunctions allowed, although this is not straightforward since it is defined using the quantifier rank. Let us see how this is the case. For example
$$
\underbrace{\bigwedge_{\xi<\delta} (\forall \vec x\uhr \xi)\, \psi_\xi(\vec x\uhr \xi)}_{ =\psi }
$$
is a sentence with $<\delta^+$ conjunction and quantifying over $<\delta^+$ variables. If the $\psi_\xi$ have quantifier rank $\leq\delta$, then such formula is an element of $\LA_{\gd^+,\gd}'$. Also, such formula is an element of $\LA_{\delta^+\delta^+}$, whenever $\psi_\xi$ is an element from $\LA_{\delta^+\delta^+}$. But, if $\{\psi_\xi:\xi<\gd \} \subset \LA_{\gd^{++}\gd^+}\setminus \LA_{\gd^+\gd^+}$  then $\psi$ cannot be an element of $\LA_{\gd^+\gd^+}$. Yet, it could occur $\psi\in\LA_{\gd^+\gd}'$.

\item Consider an infinitary vocabulary $\tau$. That is,  at least one of the extralogical symbols has infinitely many arguments. Let us say that $R_0(\vx)\in \tau$ is a relation symbol that relates $\omega$ elements, which is to say that $R^\A_0\subset A^\omega$ for any $\tau$-structure $\A$. By definition $\e\vx R_0(\vx) \in \LA_{1\, 1}'(\tau)$, but $\e\vx R_0(\vx)$ is not even an element of the logic $\LA_{\infty\go}(\tau)$. 

\item Now, if we take and $\omega_1$-ary relation $R_1$, and fix $\vfi(\vy)$ as $\e \vx R_1(\vx\vy)$ where $\ell(\vx)=\omega$, then $\vfi(\vy)\in \LA_{1\, \omega_1}'$. 

\item Even when considering a first-order vocabulary $\tau$ with at least $\omega_1$ many relations, say $R_\xi(x_0^\xi,...,x_{k_\xi}^\xi)$ $\xi<\omega_1$, we can consider
$$
\bigwedge_{\xi<\omega_1} \e x_0^\xi R_\xi(\vec x)
$$
which is not a formula of $\LA_{\omega_2\omega_1}(\tau)$ but it is an element of $\LA'_{\omega_1\omega_2}(\tau)$.
%

\end{enumerate}
\end{ej}

Whenever $\A$ and $\B$ are $\tau$-structures, we write $\A\prec'_{\gk\gl}$ for $\A\subset \B$ and for any $\vfi\in \LA'_{\gk\gl}(\tau \cup \{c_a:a \in A\})$ we have
$$
(\A,a)_{a\in A} \models \vfi \ssi (\B,a)_{a\in A}\models \vfi.
$$
Where $c_a$ is to be interpreted as $a\in A$ for both of the models $\A,\B$.

\vspace{2em}

Our goal is to obtain a generalized form of the Gap-$n$ transfer problem. First we will tackle this using the constructible hierarchy in the set-theoretic vocabulary and using only one predicate. For this section, we assume the existence of $0^\sharp$. The set $0^\sharp$ has many presentations, one of them states the existence a proper class $C$ for \emph{cardinal indiscernibles} of the class model $L$. The class $(X,<)$ is indiscernible for the (class) model $\A$ if for every increasing sequences $x_0<...<x_k$ and $y_0<...<y_k$ from $X$ we have 
$$
\A\models \vfi(x_0,...,x_k) \lra \vfi(y_0,...,y_k).
$$
One of the main features of $0^\sharp$ is that for each $\gk,\gl\in C$ we have $L_\gk \prec L_\gl$ whenever $\gk\leq \gl$. A feature that will be exploided in a stronger sense. See \cite{Silver} for more on $0^\sharp$.


\begin{thm}\label{two-case}
Assume $0^\sharp$. Let $\kappa< \lambda<\mu<\nu$ be cardinals such that  $cf(\lambda),cf(\mu),cf(\nu)\geq\kappa$. Then, 
$$
(L_\mu, L_\lambda,\in)\prec'_{\lambda\kappa} (L_\nu,L_\lambda,\in).
$$
\end{thm}
\begin{proof}

From \cite[Theorem 2]{Chang1} we can conclude the following
$$
(L_\lambda,\in)\prec'_{\lambda\kappa} (L_\mu,\in),( L_\nu,\in).
$$
Also,
$$
(L_\mu,\in)\prec'_{\mu\kappa}(L_\nu,\in).
$$
Particularly $(L_\mu,\in)\prec'_{\lambda\kappa} (L_\nu,\in)$. With this in mind, we can safely proceed to prove the theorem focusing on relativized sentences to $L_\lambda$. Take $\vfi\in \LA_{\lambda\kappa}'(\{\in,c_x:x\in L_\mu\})$ and we prove the following
$$
(L_\mu,L_\lambda,\in) \models \vfi^{L_\lambda} \ssi (L_\nu,L_\lambda,\in)\models \vfi^{L_\lambda}.
$$
But $\vfi^{L_\gl}$ is a $\gS_0$ formula and $L_\gl, L_\mu, L_\nu$ are transitive set-theoretic models. 
\end{proof}

As a corollary we obtain a particular instance of the Gap-$n$ problem.

\begin{coro}
Suppose $\gk<\gl\leq\mu<\nu$ satisfy $ cf(\mu), cf(\nu) \geq \gk$. For any $n<\omega$ we have
$$
(L_{\mu^{(n)}}, L_\mu,\in) \prec'_{\gl\gk} (L_{\nu^{(n)}},L_\nu,\in).
$$
\end{coro}

We wish to generlize this notion to include more than one predicate. In that direction we have the following theorem.

\begin{thm}\label{ncase}
Let $n<\omega$ and $\kappa,\lambda\leq\mu_0<...<\mu<\nu$ be such that $\kappa\leq cf(\lambda),cf(\mu_i),cf(\mu), cf(\nu)$ for each $i<n$. We obtain
$$
(L_{\mu_{n-1}},...,L_{\mu_0},\in)\prec'_{\kappa\lambda}  (L_{\nu_{n-1}},...,L_{\mu_0},\in).
$$
\end{thm}
\begin{proof}
We proceed by induction on $n<\go$. Suppose it holds for $n$. As in \ref{two-case}, it is enough to prove it for relativized formulas with the presentation  $\vfi^{(L_{\mu_{n-1}},...,L_{\mu_0},\in)}$. We abbreviate this to $\vfi^{\vec L}$, Notice the following list of equivalences.
$$
\mtx{
(L_{\mu_n},L_{\mu_{n-1}},...,L_{\mu_0},\in)\models \vfi^{\vec L} & \ssi \\
(L_{\mu_{n-1}},...,L_{\mu_0},\in) \models \vfi & \ssi \\
(L_{\nu_n},L_{\mu_{n-1}},...,L_{\mu_0},\in)\models \vfi^{\vec L}.
}
$$

Notice that this argument strongly depends on the fact $(L_\lambda,\in)\prec'_{\lambda\kappa} (L_\mu,\in),( L_\nu,\in)$.
\end{proof}

Thus, the principle $0^\sharp$ solves many cardinal transfer properties. For example, we obtain 
$$
(\gg^{(n)},\gg) \to_{T'} (\gd^{(n)},\gd),  
$$ 
for regular cardinals $\gd>\gg$ such that $\gg>\gl>\gk$. This includes the upwards Gap-$n$ problem in many instances. Here we write the subscript $T'$ to make explicit that these results are only for the theory $\mathrm{Diag}_{\LA'_{\gl\gk}}(L_{\mu_{n-1}},...L_{\mu_0},\in)$. Also, the results stated here and in the rest of this work are upwards. That is, we start with a ``small'' model and obtain larger models. Elsewhere, \cite{VaVi22}, a downwards version for infinitary logic with finite quantification is considered. Although this resutls are important by themselves, because of the significance of the constructible hierachy, we strive for stronger results.

\section{Stronger Gaps}

In this section we define our generalized Gap problem and we will prove that under the existence of a special first-order theory $\Phi$, we can obtain prove this generalized Gap problem in many instances.

Following Theorem \ref{ncase} we introduce the following generalization for the Gap-$n$ problem.

Let $n_0<...< n_k$ \emph{Gap-$(n_0,...,n_k)$ problem for $\LA'_{\gk\gl}$} we mean the arrow relation
\begin{equation}\label{gapnnn}
(\mu^{(n_k)}, \mu^{(n_{k-1})},...,\mu^{(n_0)}) \to_{\LA'_{\gk\gl}} (\nu^{(n_k)}, \nu^{(n_{k-1})},...,\nu^{(n_0)}).
\end{equation}

Even when  $V\neq L$, any $(\mu,\mu_k,...,\mu_0)$-model can be represented as
\begin{equation}
\A= ( L_\mu,L_{\mu_k},...,L_{\mu_0},\ldots) \tag{\ding{67}}.
\end{equation}
This is due to the fact that $|L_\ga|=|\alpha|$ for any ordinal $\alpha$.
In the following, we will assume that any $(\mu,\mu_k,...,\mu_0)$-model has this exact presentation.

We will use the relativized hierarchy $L[x]$ for some set $x$. The construction of $L[x]$ is carried out analogously to the construction of $L$, but this time one uses $x$ as a predicate to construct the levels using recursion. See \cite{NSV} for more details.
As a convention, we fix $a = \bigsqcup_{\varsigma\in \tau} \varsigma^\A$. That is, $a$ is the disjoint union of the $\A$-interpreted symbols of the vocabulary $\tau$. Suppose that $a$ satisfies condensation. This is,
\begin{center}
\emph{for any primitive recursive ordinal $\gamma$ over the predicate $a$ and \\any $M \prec_{\gS_1} L_\gamma [a]$ there is $\rho\leq \gamma$ such that $M \cong L_\rho [a]$.}
\end{center}

The theory $\Phi$ we will be dealing with is an \emph{Ehrenfeucht-Mostowski blueprint}. An EM blueprint is a theory $\Phi$ of first-order sentences that will allow us to derive (a class of) indiscernibles. In general, EM blueprints vary. In our case, we look for the following traits for $\Phi$.

\begin{enumerate}
\item $\Phi$ is a first-order theory, in some extension of the language of set theory that includes Skolem functions.
\end{enumerate}
For any ordinal $\ga\leq \omega$,
\begin{enumerate}[start =2]
\item there is $\A\models \Phi$ and $I\subset A$ a set of indiscernibles with $\suc{ I}_\A = \A$,
\item $otp(I) = \alpha$ and $|\alpha| = |A|$,
\item $\A$ is well-founded,
\item $\A\models V=L[a]$,
\item $I$ is closed and unbounded in $Or^\A$,
\item For any Skolem term $t$, we have  $t(\vx,\vec y)<y_1 \rightarrow t(\vx,\vec y) = t(\vx,\vec z)\in\Phi$. 
\item If $\omega\leq\beta<\alpha$ and $I = \{i_\xi:\xi<\alpha\}$ then $\suc{\{i_\xi:\xi<\beta\}}_\A$ is the unique model $\B\models \Phi$ which satisfies the previous items, but restricted to $\B$.
\end{enumerate}

Among the many presentations of $0^\sharp$ one can find that it is an EM blueprint. Consult \cite{Silver} for further details.

Our Gap-$(n_0,...,n_k)$ Theorem will be a consequence of the following more general theorem.

\begin{thm}\label{main}
 For any first-order vocabulary $\tau$ of power $<\delta$ and which includes the predicates $U_0,...,U_{n-1} \in \tau$ and any cardinals $\nu>\mu>\mu_{n-1}>...>\mu_0\geq\lambda>\kappa\geq\delta$, with  $\kappa\leq cf(\lambda),cf(\mu_i),cf(\mu), cf(\nu)$ for each $i<n$, we deduce
$$
(\mu,\mu_{n-1},...,\mu_0)\to_{\LA_{\lambda\kappa}'}(\nu,\mu_{n-1},...,\mu_0).
$$
\end{thm}




Before proving the theorem, we need a few facts of $L[a]$ and its EM blueprint.


\begin{af}[1]
Let $X$ is the class of indiscernibles obtained from the theory $\Phi$. Then, for any $\mu$,  $\suc{X\uhr\mu}_{L_\mu[a]} = L_\mu[a]$.
\end{af}
\begin{prueba}
This is a classical argument and is full of details. We will give a sketch that focuses on our particular scenario. Thus, it is enough to prove that any $y\in L_{\mu}[a]$ is $\{a,\in\}$-definable with $X\cap\mu$.

The main key is that the only $(\Phi,\mu)$-model $\A$ is in fact $L_\mu[a]$. Let us say that $I = \{i_\xi:\xi<\mu\}$ is the set of indiscernibles associated to $\A$. Since $a$ satisfies condensation, we know that $\A\cong L_\beta[a]$, for some ordinal $\beta$. Obviously, $\beta\geq\mu$ because $otp(I)=\mu$. If we had $\beta>\mu$, then there is $\xi<\mu$ with $\mu<i_\xi$. Therefore, the $(\Phi,\xi)$-model $\suc{\{i_\zeta:\zeta<\xi\}}_\A$ has size $\mu$, but $\xi<\mu$. Hence it has power $<\mu$. This poses a clear contradiction. Hence $\beta=\mu$.\end{prueba}

\begin{nota} Remember that we assumed that $a$ satisfies condensation. An hypothesis that might seem restrictive. But Theorem \ref{main} does not makes explicit such hypothesis. This is because it is not really restrictive for our case as stated by the following fact. If $x$ is a set, then $x$ satisfies condensation over $rk(x)$. That is, if $\ga\geq rk(x)$ and $M\prec_{\gS_1} L_\ga$ then there is $\xi\leq\ga$ such that $M\cong L_\xi[x]$. More, whenever $x$ is a set, we can find a set of ordinals $x'$ such that $L[x] = L[x']$ and $x'$ does satisfy condensation. The details can be consulted in \cite{NSV}. 
\end{nota}

With these tools and clarifications at hand, we proceed to prove Theorem \ref{main}.

\begin{proof}[Proof of \ref{main}]
We begin with the following claim. Let us remember that $\A$ has the presentation $(L_\mu,L_{\mu_{n-1}},...,L_{\mu_0},...)$.

\begin{af}[2]
For any $\vfi\in \po(\tau)$ there is some $\vfi^*\in \po(\{\in,a\}\cup\{c_\xi:\xi<\mu\})$ such that for any $\vec e\in L_\mu$.
$$
\A\models \vfi[\vec e] \ssi (L_\mu[a],\in,a)\models \vfi^*[\vec e].
$$
\end{af}
\begin{prueba}
We proceed by induction on the construction of $\vfi\in\po(\tau)$. Before that, we must deal with the terms of the vocabulary.  Since $L_\mu\subset L_\mu[a]$, the equality sing must be preserved directly. That is, for $e_0,e_1\in L_\mu$ we have $\A\models e_0=e_1$ when and only when $(L_\mu[a],\in,a)\models e_0=e_1$. Thus we define $(x=y)^*$ as $x=y$. Suppose we have an enumeration of $\tau=\{\varsigma_\xi:\xi<\rho\}$ such that $a= \bigcup_{\xi<\rho} \{(\xi,\vec x): \vec x\in \varsigma_\xi^{\mf A}\}$. Suppose  $c\in\tau$ is a constante symbol and it is the $\xi$th symbol in $\tau$. Notice the following equivalences.
$$
\mtx{
\mf A\models c=e \ssi\\
 c^\mf A =e \ssi \\
 (\xi,c^\mf A)\in a \text{ y } (\xi,c^\mf A) = (\xi,e) \ssi \\
(L_\mu[a],\in,a)\models a(\xi,c^\mf A)\wedge (\xi,c^\mf A) = (\xi,e).}
$$
Assume that $f\in \tau$ is a function symbol and it is listed as the $\zeta$th symbol of $\tau$.
$$
\mtx{
\mf A\models f(\vec e) = e_0 \ssi \\
f^\mf A(\vec e) = e_0 \ssi \\
 (\zeta,(\vec e,e_0))\in a   \ssi \\
(L_\mu[a],\in, a)\models  a\big(\zeta,(\vec e,e_0) \big).
}
$$
Now, suppose that $\vfi(\vec y)$ is $t(\vec y) = y_{n+1}$ and $t$ is a more complex term on our vocabulary. Say $t=t_{n+1}(t_0,...,t_n)$, where $t_{n+1}$ is the $\gamma$th symbol of our vocabulary and for $t_i$ there is $\vfi^*_i$ such that
$$
\mf L\models t_i(\vec e) = e_0 \ssi (L_\mu[a],\in,a)\models \vfi_i^*(\vec e,e_0).
$$
Now, define $\vfi^*$ as
$$
\exists \vec x,x_{n+1} \big( a(\gamma,(\vec x,x_{n+1})) \wedge \bigwedge_{i\leq n+1} \vfi_i^*(x_i,\vec y) \big).
$$
This takes care of the atomic formulas. Take a relation symbol $r\in\tau$, which is listed as the $\eta$th symbol in $\tau$ and $t_i$, $i<m$ are terms. Define $r^*(\vec t)$ by letting
$$
\exists \vec x \big( a(\eta,(\vec x)) \wedge \bigwedge_{i<m} \vfi^*_i(\vec x,\vec y)\big)
$$
where $\vfi^*_i(x_i,\vec y) = (x_i = t_i(\vec y))^*$. Now we define the inductive steps in a natural way:
$$
\mtx{
\neg \vfi \mapsto \neg \vfi^*, \\
\vfi\wedge\psi \mapsto \vfi^*\wedge \psi^*, \\
\exists x \vfi \mapsto \exists x \vfi^*.
}
$$
This closes the construction and we can see that $\vfi\mapsto \vfi^*$ is the required traduction.
\end{prueba}

\begin{af}[3]
If $X\subset L_{\mu}[a]$ is the canonical set of de indiscernibles $\suc{X} = L_\mu[a]$, then $X$ is also a set of indiscernibles for $\A=(L_\mu,L_{\mu_{n-1}},...,L_{\mu_0},...)$.
\end{af}
\begin{prueba}
Let $\xi_1<...<\xi_n,\zeta_1<...<\zeta_n$ be elements from $X$ and $\vfi\in\po(\tau)$. Then, consider $\vfi^*$, we know
$$
(L_\mu[a],\in,a)\models \vfi^*(\vec\xi) \lra \vfi^*(\vec\zeta).
$$
Therefore
$$
\mf A\models \vfi(\vec \xi) \lra \vfi(\vec\zeta).
$$
\end{prueba}

Now, following \cite[Propositions 2]{Chang1} 
$$
X\uhr \mu \prec_{\lambda\kappa}' X\uhr \nu.
$$
Using \cite[Proposition 1]{Chang1}, we obtain
$$
\suc{X\uhr \mu}_\A\prec'_{\lambda\kappa}\suc{X\uhr \nu}_\A.
$$
But $\suc{X\uhr \gamma}_\A = L_\gamma[a]$ for any $\delta$. At this point, we extend  $\vfi\mapsto \vfi^*$ to include the infinitary cases: $\bigwedge \Phi \mapsto \bigwedge_{\vfi\in \Phi} \vfi^*$ and $\e \vx \vfi(\vx) \mapsto \e\vx \vfi^*(\vx)$. For $\vfi\in\LA'_{\lambda\kappa}(\tau)$ and $\vec e\in L_\mu$ we obtain

$$
\mtx{ 
\mf A \models \vfi[\vec e] \ssi \\
(L_\mu[a],\in,a)\models \vfi^*[\vec e] \ssi \\
(L_\nu[a],\in,a)\models \vfi^*[\vec e] \ssi \\
(L_\nu, L_{\mu_{n-1}},...,L_{\mu_0},...)\models \vfi[\vec e].
}
$$

This is exactly what we needed.

\end{proof}

\begin{thm}[Infinitary Gap-$(n_0,...,n_k)$]
 If $\gk<\gl<\mu<\nu$ and $n_0\leq...< n_k$ then for every vocabulary $\tau$ with power $<\delta$ and such that $U_0,...,U_k\in \tau$ we have
$$
(\mu^{(n_k)}, \mu^{(n_{k-1})},...,\mu^{(n_0)}) \to_{\LA'_{\gk\gl}} (\nu^{(n_k)}, \nu^{(n_{k-1})},...,\nu^{(n_0)}).
$$

\end{thm}

If $\gk\leq \omega$ then this includes first-order logic. Notice that our results are upwards, but $\po$ can also move downwards and this shows that first order logic has the full force of the Gap-$(n_0,...,n_k)$ problem.

\begin{thm}[Gap-$(n_0,...,n_k)$]
 Suppose $n_0\leq...< n_k$ then for every vocabulary $\tau$ with power $<\delta$ and such that $U_0,...,U_k\in \tau$ we have
$$
(\mu^{(n_k)}, \mu^{(n_{k-1})},...,\mu^{(n_0)}) \to_{\po} (\nu^{(n_k)}, \nu^{(n_{k-1})},...,\nu^{(n_0)}).
$$
\end{thm}

Thus, outside of $L$, the consistency of the Gap-$n$ is positive.




\end{document}